\documentclass[11pt,a4paper]{article}

\usepackage{amsmath}
\usepackage{amsthm}
\usepackage{amssymb}
\usepackage{amscd}

\title{On Effective Non-vanishing of Weil Divisors \\
on Algebraic Surfaces
}
\author{Qihong Xie}
\date{}
\theoremstyle{plain}
\newtheorem{prop}{Proposition}
\newtheorem{lem}[prop]{Lemma}
\newtheorem{thm}[prop]{Theorem}

\newtheorem{conj}[prop]{Conjecture}
\newtheorem{prob}[prop]{Problem}

\theoremstyle{definition}
\newtheorem{defn}[prop]{Definition}
\newtheorem*{ack}{Acknowledgment}

\theoremstyle{remark}
\newtheorem{rem}[prop]{Remark}
\newtheorem{ex}[prop]{Example}

\newcommand{\Q}{\mathbb Q}
\newcommand{\R}{\mathbb R}
\newcommand{\C}{\mathbb C}

\newcommand{\A}{\mathbb A}
\newcommand{\PP}{\mathbb P}
\newcommand{\OO}{\mathcal O}
\newcommand{\RR}{\mathcal R}
\newcommand{\LL}{\mathcal L}
\newcommand{\al}{\alpha}

\begin{document}
\unitlength=1mm
\maketitle

\begin{abstract}
We give a counterexample and some conclusions for effective non-vanishing of Weil divisors on algebraic surfaces.
\end{abstract}

The basepoint-free theorem plays a crucial role in the Minimal Model Program since it yields a basepoint-free complete linear system, hence an extremal contraction morphism. More precisely, it says that if $D$ is a nef Cartier divisor such that the difference between $D$ and the log canonical divisor is nef and big, then the linear system $|mD|$ is basepoint-free for sufficiently large integer m (cf.\ \cite{kmm, km}). Obviously, it is natural and interesting to investigate the non-emptyness of the linear system $|D|$. Therefore, the following so-called Effective Non-vanishing Conjecture has been put forward formally by Yujiro Kawamata (cf.\ \cite{ka}).

\begin{conj}\label{A:1}
Let $X$ be a complete normal variety, $B$ an effective $\R$-divisor on $X$ such  that the pair $(X,B)$ is Kawamata log terminal, and $D$ a Cartier divisor on $X$. Assume that $D$ is nef and that $D-(K_X+B)$ is nef and big. Then $H^0(X,D)\neq 0$.
\end{conj}

By the Kawamata-Viehweg vanishing theorem (cf.\ \cite{kmm}), we have $H^p(X,D)\\=0$ for any positive integer $p$. Thus the condition $H^0(X,D)\neq 0$ is equivalent to saying that $\chi(X,D)\neq 0$. Since the Riemann-Roch formula of the line bundle $\OO_X(D)$ depends only on the intersections of Chern classes of $X$ and $D$, we can say that the Effective Non-vanishing Conjecture is indeed a topological problem in a sense of complex geometry. This conjecture is of great importance by at least two reasons. At first, it provides a good tool to study the geometry of Fano varieties (cf.\ \cite{ambro}). Secondly, it is closely related to the general elephant problem (cf.\ \cite{reid}).

Kawamata has proven that the Effective Non-vanishing Conjecture holds for all log surfaces with only Kawamata log terminal singularities (cf.\ \cite{ka}). In dimension three, this conjecture is more complicated, only a few results are known.

It is also interesting to consider a new problem by replacing Cartier with $\Q$-Cartier.

\begin{prob}\label{A:2}
Notation are the same as in Conjecture \ref{A:1}. Assume that $D$ is a nef $\Q$-Cartier Weil divisor on $X$, such that $D-(K_X+B)$ is nef and big. Does $H^0(X,D)\neq 0$ hold?
\end{prob}

Unfortunately, this problem does not hold generally. In dimension three, Fletcher has given the following famous example (cf.\ \cite{fle}).

\begin{ex}\label{A:3}
A general weighted complete intersection $X_{12,14}$ in $\PP(2,3,4,5,6,7)$ is a terminal $\Q$-Fano threefold of Picard number one with $H^0(-K_X)=0$. This variety has the following isolated singularities: one of type $\frac{1}{5}(4,1,2)$, two of type $\frac{1}{3}(2,1,1)$, and seven of type $\frac{1}{2}(1,1,1)$.
\end{ex}

In dimension two, such a counterexample of Problem \ref{A:2} also exists, which is Example \ref{A:7} given later in this paper. At first, we give some notation and definitions for the convenience of the reader.

\begin{defn}\label{A:4}
A normal surface $X$ is said to be log terminal (resp.\ canonical), if the following conditions are satisfied.

(1) $K_X$ is a $\Q$-Cartier divisor, i.e., there exists some positive integer $r$ such that $rK_X$ is Cartier, and

(2) For any resolution $f:Y\rightarrow X$, if we write $K_Y=f^*K_X+\sum\al_iE_i$, where $E_i$ are irreducible exceptional curves of $f$, then $\al_i>-1$ (resp.\ $\al_i\geq 0$) for all $i$.
\end{defn}

The index of $X$ is, by definition, the smallest positive integer $r$ such that $rK_X$ is Cartier. If the index of $X$ is one, then $X$ is said to be Gorenstein. If the resolution $f$ is minimal, i.e., $K_Y$ is $f$-nef, then we have that $\al_i\leq 0$ for all $i$ since the intersection matrix $(E_i.E_j)$ is negative definite.

\begin{defn}\label{A:5}
A normal surface $X$ is called a log del Pezzo surface, if $X$ has only log terminal singularities and $-K_X$ is ample.
\end{defn}

For a $\Q$-divisor $D=\sum a_iD_i$, the round up (resp.\ round down) of $D$ is defined to be $\lceil D\rceil=\sum\lceil a_i\rceil D_i$ (resp.\ $\lfloor D\rfloor=\sum\lfloor a_i\rfloor D_i$), where $\lceil a\rceil$ (resp.\ $\lfloor a\rfloor$) is the smallest integer not less (resp.\ the largest integer not greater) than $a$.

\begin{prop}\label{A:6}
Let $X$ be a log terminal surface, $f:Y\rightarrow X$ the minimal resolution of $X$. Assume that $K_Y=f^*K_X+\sum\al_iE_i$, where $-1<\al_i\leq 0$ for all $i$. Then for each integer $n\geq 0$, $\chi(-nK_X):=\chi(X,\OO_X(-nK_X))$ is equal to
\[ {n(n+1)\over 2}K_Y^2-{2n+1\over 2}K_Y.(\sum\lceil(n+1)\al_i\rceil E_i)+{1\over 2}(\sum\lceil(n+1)\al_i\rceil E_i)^2+\chi(\OO_Y). \]
\end{prop}

\begin{proof}
Let $D=-nK_X$. Then we have
\begin{eqnarray}
f^*(D-K_X) & = & -(n+1)f^*K_X=-(n+1)(K_Y-\sum\al_iE_i) \nonumber \\
           & = & (-nK_Y+(n+1)\sum\al_iE_i)-K_Y. \nonumber
\end{eqnarray}

We define $\widetilde{D},\widetilde{B}$ as follows:
\begin{eqnarray}
\widetilde{D} & = & \lceil-nK_Y+(n+1)\sum\al_iE_i\rceil=-nK_Y+\sum\lceil(n+1)\al_i\rceil E_i, \nonumber \\
\widetilde{B} & = & \sum(\lceil(n+1)\al_i\rceil-(n+1)\al_i)E_i. \nonumber
\end{eqnarray}

It is easy to see that $f^*(D-K_X)=\widetilde{D}-(K_Y+\widetilde{B})$ is an $f$-nef and $f$-big divisor, and $\widetilde{B}$ has simple normal crossing support. By the Kawamata-Viehweg vanishing theorem, we have $R^jf_*\OO_Y(\lceil\widetilde{D}-\widetilde{B}\rceil)=R^jf_*\OO_Y(\widetilde{D})=0$ for any positive integer $j$.

Since we have
\begin{eqnarray}
f^*D=-nK_Y+\sum n\al_iE_i=\widetilde{D}+\sum(n\al_i-\lceil(n+1)\al_i\rceil)E_i,
\nonumber \\
n\al_i-\lceil(n+1)\al_i\rceil < (n+1)\al_i-\lceil(n+1)\al_i\rceil+1\leq 1, \nonumber
\end{eqnarray}
there exists a $\Q$-divisor $H=\sum h_iE_i$ on $Y$ such that $0\leq h_i<1$ and $h_i-(n\al_i-\lceil(n+1)\al_i\rceil) \geq 0$ for each $i$. For a $\Q$-divisor $M$ on $Y$, the sheaf $\OO_Y(M)$ is defined to be $\OO_Y(\lfloor M\rfloor)$. Thus $f_*\OO_Y(\widetilde{D})=f_*\OO_Y(\widetilde{D}+H)=\OO_X(D)$ by the projection formula. It follows from \cite{bre} that we have the Leray spectral sequence
\[ E^{i,j}_2=H^i(X, R^jf_*\OO_Y(\widetilde{D}))\Rightarrow H^{i+j}(Y, \OO_Y(\widetilde{D})). \]
Since $E^{i,j}_2=0$ for all $i\geq 0$ and $j>0$, $E^{i,0}_2\cong E^{i,0}_\infty$ implies that $h^i(X,\OO_X(D))=h^i(Y,\OO_Y(\widetilde{D}))$ for all $i\geq 0$, hence $\chi(X,\OO_X(D))=\chi(Y,\OO_Y(\widetilde{D}))$. Therefore we have
\begin{eqnarray}
&   & \chi(-nK_X) \nonumber \\
& = & \chi(X,\OO_X(D))=\chi(Y,\OO_Y(\widetilde{D}))={1\over 2}\widetilde{D}(\widetilde{D}-K_Y)+\chi(\OO_Y) \nonumber \\
& = & {1\over 2}(-nK_Y+\sum\lceil(n+1)\al_i\rceil E_i)(-(n+1)K_Y+\sum\lceil(n+1)\al_i\rceil E_i)+\chi(\OO_Y) \nonumber \\
& = & {n(n+1)\over 2}K_Y^2-{2n+1\over 2}K_Y.(\sum\lceil(n+1)\al_i\rceil E_i)+{1\over 2}(\sum\lceil(n+1)\al_i\rceil E_i)^2+\chi(\OO_Y). \nonumber
\end{eqnarray}
\end{proof}

\begin{ex}\label{A:7}
There exists a log del Pezzo surface $S$ of Picard number one such that $H^0(S,-K_S)=0$. 

Let $B$ be a smooth conic in $\PP^2$, $D$ a tangent line to $B$ at $d$, and $A$ a secant line intersecting $B$ at $\{a,b\}$. Let $g_1: S_1\rightarrow \PP^2$ be the blow-up of $\PP^2$ along $d\in D$ with an exceptional curve $L_1$. By abuse of notation, we always denote the strict transform of a curve by the same symbol. Let $g_i: S_{i+1}\rightarrow S_i$ be the blow-up of $S_i$ along $d\in D$ with exceptional curves $L_{i+1}$ for $i=1,2$. It is easy to see that $L_1^2=L_2^2=-2, L_3^2=-1$. In a similar way, let $g_4: S_4\rightarrow S_3$ be the blow-up of $S_3$ 5 times along $b\in B$. Assume that the exceptional curves of $g_4$ are denoted by $M_1,M_2,M_3,M_4,M_5$ respectively with $M_1^2=M_2^2=M_3^2=M_4^2=-2$, and $M_5^2=-1$. Let $g_5: \widetilde{S}\rightarrow S_4$ be the blow-up of $S_4$ 5 times along $a\in A$. Assume that the exceptional curves of $g_5$ are denoted by $N_1,N_2,N_3,N_4,N_5$ respectively with $N_1^2=N_2^2=N_3^2=N_4^2=-2$, and $N_5^2=-1$. Denote by $g: \widetilde{S}\rightarrow\PP^2$ the composition of all $g_i$.

We consider two reduced connected curves $X_1:=D\cup A\cup M_1\cup M_2\cup M_3\cup M_4$ and $X_2:=L_1\cup L_2\cup B\cup N_1\cup N_2\cup N_3\cup N_4$ on $\widetilde{S}$. It is easy to verify that the intersection matrix associated with $X_i$ is negative definite, and the arithmetic genus of any positive cycle with support on $X_i$ is non-positive for $i=1,2$. Thus $X_1\cup X_2$ is contractible, and there is a contraction $f: \widetilde{S}\rightarrow S$ such that $S$ is a normal projective surface of Picard number one with only two singular points $x_1,x_2$, which is the contraction point of $X_1,X_2$ respectively (cf.\ \cite{art}).

It is easy to see that $\widetilde{S}$ is just the minimal resolution of $S$, and the dual graph $\Gamma$ consists of two connected components corresponding to $x_1,x_2$ respectively.

\noindent $\Gamma$:\quad
\begin{picture}(50,10)(0,0)
\put(0,0){\circle*{1}}
\put(-1,2){$2^D$}
\put(0,0){\line(1,0){10}}
\put(10,0){\circle*{1}}
\put(9,2){$5^A$}
\put(10,0){\line(1,0){10}}
\put(20,0){\circle*{1}}
\put(19,2){$2^{M_1}$}
\put(20,0){\line(1,0){10}}
\put(30,0){\circle*{1}}
\put(29,2){$2^{M_2}$}
\put(30,0){\line(1,0){10}}
\put(40,0){\circle*{1}}
\put(39,2){$2^{M_3}$}
\put(40,0){\line(1,0){10}}
\put(50,0){\circle*{1}}
\put(49,2){$2^{M_4}$}
\end{picture}
\qquad
\begin{picture}(60,10)(0,0)
\put(0,0){\circle*{1}}
\put(-1,2){$2^{L_1}$}
\put(0,0){\line(1,0){10}}
\put(10,0){\circle*{1}}
\put(9,2){$2^{L_2}$}
\put(10,0){\line(1,0){10}}
\put(20,0){\circle*{1}}
\put(19,2){$4^B$}
\put(20,0){\line(1,0){10}}
\put(30,0){\circle*{1}}
\put(29,2){$2^{N_1}$}
\put(30,0){\line(1,0){10}}
\put(40,0){\circle*{1}}
\put(39,2){$2^{N_2}$}
\put(40,0){\line(1,0){10}}
\put(50,0){\circle*{1}}
\put(49,2){$2^{N_3}$}
\put(50,0){\line(1,0){10}}
\put(60,0){\circle*{1}}
\put(59,2){$2^{N_4}$}
\end{picture}
\vskip 5mm
\noindent where $2^D$ denotes the $(-2)$-curve on $\widetilde{S}$ which is the strict transform of $D$, and so on.

By an elementary calculation, we have
\begin{eqnarray}
K_{\widetilde{S}} & = & f^*K_S-\frac{15}{37} D-\frac{30}{37} A-\frac{24}{37} M_1-\frac{18}{37} M_2-\frac{12}{37} M_3-\frac{6}{37} M_4 \nonumber \\
            & - & \frac{5}{19} L_1-\frac{10}{19} L_2-\frac{15}{19} B-\frac{12}{19} N_1-\frac{9}{19} N_2-\frac{6}{19} N_3-\frac{3}{19} N_4 \nonumber
\end{eqnarray}

It follows easily that $K_S^2=8/(37\times 19)>0$, then $-K_S$ is ample since $S$ is rational and $\rho(S)=1$. Furthermore, we have $H^0(S,-nK_S)=0$ for $1\leq n\leq 5$, but $H^0(S,-6K_S)\cong\C$ by Proposition \ref{A:6}.
\end{ex}

\begin{rem}\label{A:8}
For a log del Pezzo surface $X$, we define $\tau(X)$ to be the smallest positive integer $\tau$ such that $H^0(X,-\tau K_X)\neq 0$. It is easy to see that the index $r(X)$ of $X$ is not less than $\tau(X)$. The study of complements on surfaces is helpful to solve the non-vanishing of $H^0(X,-nK_X)$ for log del Pezzo surfaces. Shokurov proved that there exists a positive integer $N$ such that $\tau(X)\leq N$ for all log del Pezzo surfaces (cf.\ \cite{sho}). In particular, Prokhorov used the same method to yield that $H^0(X,-nK_X)\neq 0$ for some $n\in\RR_2:=\{1,2,3,4,6\}$, if $K_X^2>4$ (cf.\ \cite{pro}).
\end{rem}

Example \ref{A:7} illustrates that Problem \ref{A:2} does not hold generally for surfaces. But, when some reasonable condition is assumed, we hope that the effective non-vanishing of Weil divisors holds. We give the following theorem as a special example.

\begin{thm}\label{A:9}
Let X be a Gorenstein log del Pezzo surface of Picard number one, $D$ a nef and big $\Q$-Cartier Weil divisor on X. Then $H^0(X,D)\neq 0$.
\end{thm}

It is easy to see that $X$ has only canonical singularities. The condition $\rho(X)=1$ enables us to list all types of its singularities.

\begin{lem}\label{A:10}
Let X be a Gorenstein log del Pezzo surface with $\rho(X)=1$. Then its singularity type is one of the following:

\noindent $A_1$, $A_1+A_2$, $A_4$, $2A_1+A_3$, $D_5$, $A_1+A_5$, $3A_2$, $E_6$, $3A_1+D_4$, $A_7$, $A_1+D_6$, $E_7$, $A_1+2A_3$, $A_2+A_5$, $D_8$, $2A_1+D_6$, $E_8$, $A_1+E_7$, $A_1+A_7$, $2A_4$, $A_8$, $A_1+A_2+A_5$, $A_2+E_6$, $A_3+D_5$, $4A_2$, $2A_1+2A_3$, $2D_4$.
\end{lem}

\begin{proof}
$X$ is said to be of type, for example, $2A_1+A_3$, if $X$ has three singular points, two of type $A_1$, and one of type $A_3$. We refer to \cite{mz} for the proof.
\end{proof}

\begin{defn}\label{A:11}
The quotient $X=\A^n/\mu_r$ is said to be of cyclic quotient singularity of type $\frac{1}{r}(a_1,\cdots,a_n)$, if $\mu_r$ acts on $\A^n$ by
\[ \mu_r\ni\varepsilon:(x_1,\cdots,x_n)\mapsto(\varepsilon^{a_1}x_1,\cdots,\varepsilon^{a_n}x_n). \]
Let $\pi: \A^n\rightarrow X$ be the quotient morphism. Then the group $\mu_r$ acts on $\pi_*\OO_{\A^n}$, and so decomposes it into $r$ eigensheaves $\LL_i=\{f|\varepsilon(f)=\varepsilon^i\cdot f$ for all $\varepsilon\in\mu_r \}$ for $i=0,\cdots, r-1$. A singularity $p\in X$ with a Weil divisor $D$ is said to be of cyclic quotient singularity of type $_i(\frac{1}{r}(a_1,\cdots,a_n))$, if locally, $p\in X$ is isomorphic to a point of type $\frac{1}{r}(a_1,\cdots,a_n)$, and $\OO_X(D)\cong\LL_i$.
\end{defn}

There are two keypoints in the proof of Theorem \ref{A:9}. First, one canonical singularity can be considered as the cyclic quotient of another canonical singularity by using the local cyclic cover. Second, we can calculate the Euler characteristic of a Weil divisor, thanks to the Singular Riemann-Roch Formula.

The following lemma is helpful, which can be found in \cite{reid}.

\begin{lem}\label{A:12}
Let $Q\in Y:(f=0)\subset \A^3$ be a surface germ. Assume that there is an action of $\mu_r$ on $Y$ free outside $Q$. Let $p\in X = (Q\in Y)/\mu_r$ be the surface germ and $\pi :Y\to X$ the corresponding cyclic cover. If $p\in X$ is canonical, then $Q\in Y$ is also canonical and $\pi :Y\to X$ is exactly one of the following types.

\begin{center}
\begin{tabular}{|c|c|l|l|l|} \hline
    & r & Type & f & Description \\ \hline
{\rm (1)} & {\rm any} & ${1\over r}(1,-1,0)$ & $xy+z^n$ & $A_{n-1}\stackrel{r:1}{\to}A_{rn-1}$ \\
{\rm (2)} & {\rm 4} & ${1\over 4}(1,3,2)$ & $x^2+y^2+z^{2n-1}$ & $A_{2n-2}\stackrel{4:1}{\to}D_{2n+1}$ \\ 
{\rm (3)} & {\rm 2} & ${1\over 2}(0,1,1)$ & $x^2+y^2+z^{2n}$ & $A_{2n-1}\stackrel{2:1}{\to}D_{n+2}$ \\ 
{\rm (4)} & {\rm 3} & ${1\over 3}(0,1,2)$ & $x^2+y^3+z^3$ & $D_4\stackrel{3:1}{\to}E_6$ \\ 
{\rm (5)} & {\rm 2} & ${1\over 2}(1,1,0)$ & $x^2+y^2z+z^{n}$ & $D_{n+1}\stackrel{2:1}{\to}D_{2n}$ \\ 
{\rm (6)} & {\rm 2} & ${1\over 2}(1,0,1)$ & $x^2+y^3+z^4$ & $E_6\stackrel{2:1}{\to}E_7$ \\ \hline
\end{tabular}
\end{center}
\end{lem}

\begin{proof}[Proof of Theorem \ref{A:9}]

If $D$ is a Cartier divisor, then the non-vanishing follows from \cite{ka}. Now assume that $D$ is not Cartier. For the reflexive sheaf $\OO_X(D)$, we have the following Singular Riemann-Roch Formula (cf.\ \cite{reid}):
\[ \chi(X,\OO_X(D))={1\over 2}D.(D-K_X)+\chi(\OO_X)+\sum_{p\in X}c_p(D), \]
where $c_p(D)\in\Q$ is a contribution due to the singularity of $\OO_X(D)$ at $p$, the sum is taken over the points $p\in X$ at which $D$ is not Cartier.

If $p\in X$ and $D$ is of type $_i(\frac{1}{r}(1,-1))$, then $c_p(D)=-i(r-i)/2r$. For other cases, we use the suitable $Q$-smoothing of $X$, then there is a flat deformation $\{X_\lambda,D_\lambda\}$ of $\{X,D\}$ such that $X_\lambda$ has only cyclic quotient singularities. Thus $c_p(D)$ is equal to the sum of the contributions of cyclic quotient singularities of $X_\lambda$. The concrete conclusions are listed below.

\begin{center}
\begin{tabular}{|c|c|c|} \hline
    & $X_\lambda$ & $Q$-smoothing points \\ \hline
(1) & $f+\lambda z$ & n \\
(2) & $f+\lambda z$ & 2n+1 \\
(3) & $f+\lambda x$ & 2 \\
(4) & $f+\lambda x$ & 2 \\
(5) & $f+\lambda z$ & n \\
(6) & $f+\lambda y$ & 3 \\ \hline
\end{tabular}
\end{center}

Since $-K_X$ is an ample Cartier divisor and $D$ is a nef and big Weil divisor, it is easy to show that $\chi(\OO_X)=1$ and $D.(-K_X)\geq 1$. Thus $h^0(X,D)\geq D^2/2+3/2+\sum_{p\in X}c_p(D)$.

Next, we will estimate $c_p(D)$ for a given germ $p\in X$ at which $D$ is not Cartier. Since the cyclic covers in Lemma \ref{A:12} can be constructed locally, and the contribution $c_p(D)$ depends only on the local analytic type of $p\in X$ and $D$, we can take $\pi :Y\to X$ to be a suitable cyclic cover in Lemma \ref{A:12} and estimate $c_p(D)$ as follows:

Let $p\in X$ be of type $A_{n-1}$. Then $p\in X$ is a cyclic quotient singularity of type $\frac{1}{n}(1,-1)$. By definition, we have $c_p(D)=-i_0(n-i_0)/2n$ for some integer $i_0\in [0,n]$. If $n=2k$, then $c_p(D)\geq -k/4$. If $n=2k+1$, then $c_p(D)\geq -k(k+1)/2(2k+1)$. 

Let $p\in X$ be of type $D_{n+2}$. We may consider $p\in X$ as a cyclic quotient of type (3). Then $c_p(D)\geq -1/4\times 2=-1/2$.

Let $p\in X$ be of type $E_6$. We may consider $p\in X$ as a cyclic quotient of type (4). Then $c_p(D)\geq -1/3\times 2=-2/3$.

Let $p\in X$ be of type $E_7$. We may consider $p\in X$ as a cyclic quotient of type (6). Then $c_p(D)\geq -1/4\times 3=-3/4$.

Since $E_8$ does not appear in the table of Lemma \ref{A:12}, $p\in X$ could not be of type $E_8$ provided that $D$ is not Cartier at $p$.

It is easy to check for all cases in Lemma \ref{A:10} that $\sum_{p\in X}c_p(D)\geq -3/2$. Hence $h^0(X,D)\geq D^2/2>0$.
\end{proof}

\begin{ack}
I would like to express my gratitude to Professor Yujiro Kawamata for his valuable advice and warm encouragement. I also thank Professor S. A. Kudryavtsev for informing me some known results in Remark \ref{A:8}, and Dr.\ Hokuto Uehara for stimulating discussions on this paper.
\end{ack}

\textsc{Graduate School of Mathematical Sciences, University of Tokyo, 3-8-1 Komaba, Meguro, Tokyo 153-8914, Japan}

\textit{E-mail address}: \texttt{xqh@ms.u-tokyo.ac.jp}

\end{document}